\documentclass[a4paper]{amsart}
\usepackage{mathrsfs}
\usepackage{bbm}
\usepackage{latexsym, amssymb,amsmath,amsthm}
\usepackage[all, knot]{xy}
\xyoption{arc}
\newtheorem{thm}{Theorem}[section]
\newtheorem{cor}[thm]{Corollary}
\newtheorem{lem}[thm]{Lemma}
\newtheorem{prop}[thm]{Proposition}
\newtheorem{exa}[thm]{Example}
\newtheorem{con}[thm]{Convention}

\newtheorem{rem}[thm]{Remark}
   \def\op{\oplus} \def\ot{\otimes}
\def\Hom{\operatorname {Hom}}

\def\Ext{\operatorname {Ext}} 
 \def\k{\mathbbm{k}}

\begin{document}
\title[Ore extensions of graded 2-CY algebras]{\bf Graded 3-Calabi-Yau algebras as Ore extensions of 2-Calabi-Yau algebras}

\author{Ji-Wei He, Fred Van Oystaeyen and Yinhuo Zhang}
\address{J.-W. He\newline \indent Department of Mathematics, Shaoxing College of Arts and Sciences, Shaoxing Zhejiang 312000,
China\newline \indent Department of Mathematics and Computer
Science, University of Antwerp, Middelheimlaan 1, B-2020 Antwerp,
Belgium} \email{jwhe@usx.edu.cn}
\address{F. Van Oystaeyen\newline\indent Department of Mathematics and Computer
Science, University of Antwerp, Middelheimlaan 1, B-2020 Antwerp,
Belgium} \email{fred.vanoystaeyen@ua.ac.be}
\address{Y. Zhang\newline
\indent Department WNI, University of Hasselt, Universitaire Campus,
3590 Diepenbeek, Belgium} \email{yinhuo.zhang@uhasselt.be}

%\thanks{}
\date{}
\begin{abstract} We study a class of graded algebras obtained from Ore extensions of graded Calabi-Yau algebras of dimension 2. It is proved that these algebras are graded Calabi-Yau and graded coherent. The superpotentials associated to these graded Calabi-Yau algebras are also constructed.
\end{abstract}


\maketitle
\vspace{5mm}

\section*{Introduction}

Recently, Smith studied in \cite{Sm} a remarkable graded Calabi-Yau algebra $B$ of dimension 3 constructed from the octonions. Amongst other things, Smith proved that $B$ is a graded Ore extension of an Artin-Schelter regular algebra of global dimension 2 and uses that fact to show that $B$ is graded 3-Calabi-Yau and graded coherent.

In this note, we show that the Calabi-Yau property and the coherence of $B$ do not occur incidentally. A large class of graded algebras that are Ore extensions of graded Calabi-Yau algebras are themselves graded Calabi-Yau. The main result of this note is the following.

\begin{thm}\label{thm} Let $V$ be a finite dimension vector space with basis $\{x_1,\dots,x_n\}$, let $M$ be an invertible $n\times n$ anti-symmetric matrix, and define
 $$r=(x_1,\dots,x_n)M\left(
                                                                          \begin{array}{c}
                                                                            x_1 \\                                                                          \vdots \\
x_n                                                                         \end{array}
                                                                        \right)
\in T(V).$$ Let $A=T(V)/\langle r\rangle$, where $\langle r\rangle$ is the ideal of $T(V)$ generated by $r$. Let $\delta$ be a degree-one graded derivation of $T(V)$ such that $\delta(r)=0$. Then $\delta$ induces a graded derivation $\bar{\delta}$ on $A$. Let $B=A[z;\bar{\delta}]$ be the Ore extension of $A$ defined by $\bar{\delta}$. Then the following hold:
\begin{itemize}
\item[(i)] $B$ is a graded 3-Calabi-Yau algebra.

\item[(ii)] Let $\widehat{V}=V\op\k z$, and $Q=\left(
                                   \begin{array}{ccc}
                                     -1 & 0  \\
                                     0 & M \\
                                   \end{array}
                                 \right)$. Let        $w=(z,x_1,\dots,x_n)Q\left(
                                                                          \begin{array}{c}
                                                                            r \\
                                                                            r_1 \\
                                                                            \vdots \\
                                                                            r_n
                                                                          \end{array}
                                                                        \right)$, where $r_i=z\ot x_i-x_i\ot z-\delta(x_i)\in \widehat{V}\ot \widehat{V}$ for all $i=1,\dots,n$. Then
$(\alpha\ot1\ot1)(w)=(1\ot1\ot\alpha)(w)$ for all $\alpha\in (\widehat{V})^*$, and $A[z;\bar{\delta}]\cong T(\widehat{V})/\langle \partial_{x_i}(w):i=0,\dots,n\rangle$, where we set $x_0=z$ and $\partial_{x_i}(w)$ is the cyclic partial derivative of $w$ with respect to $x_i$.

\item[(iii)] Write $\delta(x_i)=\sum_{s,t=1}^nk^i_{st}x_i\ot x_j$ for all $i=1,\dots,n$. Assume there is an integer $j$ such that $k^i_{jj}=0$ for all $i=1,\dots,n$, and $M$ is a standard anti-symmetric matrix. Then $B$ is graded coherent.
\end{itemize}
\end{thm}

Most of this note is devoted to the proof of Theorem \ref{thm}. However, we will go a bit further to discuss the properties of the algebra $B$. Smith's algebra in \cite{Sm} is an example satisfying the conditions in the theorem. We will provide a few more examples. We remark that any quadratic algebra $A$ defined by an invertible anti-symmetric matrix as in the above theorem is isomorphic to a quadratic algebra defined by a standard anti-symmetric matrix (see Convention \ref{con} for the definition), because, for every invertible anti-symmetric matrix $M$, there is an invertible matrix $P$ such that $P^tMP$ is a standard anti-symmetric matrix.

\begin{rem}\label{rrem}{\rm Let $V$ be a finite dimensional vector space with basis $\{x_1,\dots,x_n\}$. Take an element $r\in V\ot V$. Since $V\ot V\cong \Hom_\k(V^*,V)$, the element $r$ corresponds to a linear map $f_r:V^*\to V$. The {\it rank} of $r$, denoted by rank($r$), is defined to be the rank of $f_r$ (cf. \cite[Introduction]{Z2}). One sees $$\text{rank}(r)=\min\{m|r=u_1\ot v_1+\cdots+u_m\ot v_m, \text{ for some }u_i,v_i\in V\}.$$ It has been shown that certain features of the algebra $T(V)/\langle r\rangle$ entirely depend on rank($r$) (cf. \cite[Theorem 0.1]{Z2}). If $M$ is an $n\times  n$ matrix and $r=(x_1,\dots,x_n)M\left(
                                                                          \begin{array}{c}
                                                                            x_1 \\
                                                                            \vdots \\
                                                                            x_n
                                                                          \end{array}
                                                                        \right)
\in V\ot V$, then rank($r$)=rank($M$). Therefore, the condition that $M$ is invertible in Theorem \ref{thm} is equivalent to the condition that rank($r$)=$n$. }
\end{rem}

Throughout $\k$ is a fixed field. The unadorned $\ot$ means $\ot_\k$. Let $U=\op_{n\in \mathbb{Z}}U_n$ be a graded vector space, and $l$ an integer. We write $U(l)$ for the graded vector space with degree $k$ component $U(l)_k=U_{k+l}$.

A connected graded algebra $A$ is called a {\it graded Calabi-Yau algebra} of dimension $d$, or simply {\it graded d-CY algebra} (cf. \cite{Gin}), if

\begin{itemize}
\item[(i)] $A$ is
homologically smooth; that is, $A$ has a finite resolution by
finitely generated graded projective left $A^e$-modules, where $A^e=A\ot A^{op}$ is the enveloping algebra of $A$;

\item[(ii)] the projective dimension of $A$ as a left $A^e$-module is $d$, and $\Ext^i_{A^e}(A,A\ot A)=0$ if $i\neq d$ and $\Ext_{A^e}^d(A,A\ot A)\cong A(l)$ for some integer $l$ as a right $A^e$-module.
\end{itemize}

We refer to \cite{Z} (also, cf. \cite{Ber} and \cite{DV}) for the basic properties of a graded 2-CY algebra.

\section{Ore extensions of graded Calabi-Yau algebras of dimension 2}

Let $V$ be a vector space with basis $\{x_1,\dots,x_n\}$. Let $A$ be a graded quotient algebra of $T(V)$. If $A$ is a graded 2-CY algebra, then it is defined by an $n\times n$ invertible anti-symmetric matrix $M$ \cite{Z} (also, cf. \cite[Proposition 3.4]{Ber}); that is, $A\cong T(V)/\langle r\rangle$ with $r=(x_1,\dots,x_n)M(x_1,\dots,x_n)^t$. Henceforth, we assume $A=T(V)/\langle r\rangle$ with \\ $r=(x_1,\dots,x_n)M(x_1,\dots,x_n)^t$ for some fixed anti-symmetric matrix $M$. Let $\pi:T(V)\to A$ be the natural projection map. Since degree($r$)=2, we can, and we will, identify $V$ with $A_1$ through the projection $\pi$.

Let $\delta:V\to V\ot V$ be a linear map. Then $\delta$ extends in  a unique way to a degree-one derivation (also denoted by $\delta$) of $T(V)$. If $\delta(r)\in \langle r\rangle$, then $\delta$ induces a derivation $\bar{\delta}$ on $A$.

From now on, we assume that $\delta(r)\in\langle r\rangle$. Let $B=A[z;\bar{\delta}]$ be the graded Ore extension of $A$ by the derivation $\bar{\delta}$; that is, we view $z$ as an element of degree 1, and $za=az+\bar{\delta}(a)$ for all $a\in A$.

Zhang proved in \cite{Z} that $A$ is a Koszul algebra of global dimension 2, and the minimal projective resolution of ${}_A\k$ can be written as follows:
\begin{equation}\label{eq1}
   0\longrightarrow A\ot \k r\overset{\overline{d}^{-2}}\longrightarrow A\ot V\overset{\overline{d}^{-1}}\longrightarrow A\overset{\varepsilon}\longrightarrow{}_A\k\longrightarrow0,
\end{equation}
where $\varepsilon$ is the augmentation map, $\overline{d}^{-1}(1\ot x)=\pi(x)$ for all $x\in V$, and $\overline{d}^{-2}(1\ot r)=r\in A_1\ot V$. Since $B$ is an Ore extension of $A$, $B$ is a Koszul algebra of global dimension 3.
Note that $B_A$ is a free $A$-module. Applying $B\ot_A-$ to the sequence (\ref{eq1}), we obtain the exact sequence
\begin{equation}\label{eq2}
    0\longrightarrow B\ot\k r\overset{d^{-2}}\longrightarrow B\ot V\overset{d^{-1}}\longrightarrow B\longrightarrow B/BA_{\ge1}\longrightarrow0,
\end{equation} where the unlabeled map is the natural projection map, $d^{-1}(1\ot x)=\pi(x)\in B_1$ for all $x\in V$, and $d^{-2}(1\ot r)=r\in B_1\ot V$.

\begin{lem}\label{lemr1} Suppose that $\delta(r)=0$ and let $B=A[z;\bar{\delta}]$ be as above. We have the following morphism of cochain complexes:
$$\xymatrix{B\ot \k r\ar[r]^{d^{-2}}\ar[d]_{f^{-2}} &B\ot V\ar[d]_{f^{-1}}\ar[r]^{d^{-1}}&B\ar[d]^{f^{0}}\\
 B\ot \k r\ar[r]^{d^{-2}} & B\ot V\ar[r]^{d^{-1}}& B,
}$$ where the vertical arrows are left $B$-module morphisms $f^{-2}(1\ot r)=z\ot r$, $f^{-1}(1\ot x)=z\ot x-\delta(x)$ for all $x\in V$, and $f^{0}(1)=z$.
\end{lem}
\proof We write $r=\sum_{i=1}^nu_i\ot x_i$ with all $u_i\in V$, and assume $\delta(x_i)=\sum_{j=1}^n y_{ij}\ot x_j$ for all $i=1,\dots,n$ with all $y_{ij}\in V$. We prove the commutativity of the left square. The commutativity of the right one is easy. The identity $\delta(r)=0$ is equivalent to $\sum_{i=1}^n\delta(u_i)\ot x_i+\sum_{i=1}^n u_i\ot \delta(x_i)=0$, which is in turn equivalent to $\sum_{i=1}^n\delta(u_i)\ot x_i+\sum_{i=1}^n \sum_{j=1}^nu_i\ot y_{ij}\ot x_j=0$. Applying the map $\pi\ot 1:T(V)\ot V\to A\ot V$ to the last identity, we obtain $\sum_{i=1}^n\bar{\delta}(u_i)\ot x_i+\sum_{i=1}^n \sum_{j=1}^nu_iy_{ij}\ot x_j=0$. Hence
\begin{equation}\label{eqq1}
    \bar{\delta}(u_i)=-\sum_{j=1}^nu_jy_{ji}
\end{equation}
for all $i=1,\dots,n$. The following equations hold:
\begin{eqnarray}
 \nonumber f^{-1}\circ d^{-2}(1\ot r)&=&f^{-1}(\sum_{i=1}^nu_i\ot x_i) \\
  \nonumber &=& \sum_{i=1}^nu_iz\ot x_i- \sum_{i=1}^n\sum_{j=1}^n u_i y_{ij}\ot x_j\\
 \nonumber &=& \sum_{i=1}^n(u_i z- \sum_{j=1}^n u_j y_{ji})\ot x_i,
\end{eqnarray}
and
\begin{eqnarray}
  \nonumber d^{-2}\circ f^{-2}(1\ot r) &=& d^{-2}(z\ot r)= \sum_{i=1}^nzu_i\ot x_i\\
  \nonumber &=& \sum_{i=1}^n(u_i z+\bar{\delta}(u_i))\ot x_i.
\end{eqnarray} By Equation (\ref{eqq1}), $f^{-1}\circ d^{-2}(1\ot r)=d^{-2}\circ f^{-2}(1\ot r)$. Hence the left square of the diagram commutes.
\qed

The mapping cone of the morphism in Lemma \ref{lemr1} provides a graded projective resolution of the trivial module ${}_B\k$ (see also, \cite{GS,P}).

\begin{lem} Let $r$ and $B$ be the same as in Lemma \ref{lemr1}. The minimal projective resolution of ${}_B\k$ is as follows:
$$0\longrightarrow B\ot \k r\overset{\partial^{-3}}\longrightarrow B\ot \k r\op B\ot V\overset{\partial^{-2}}\longrightarrow B\ot V\op B\overset{\partial^{-1}}\longrightarrow B\longrightarrow\k\longrightarrow0,$$ where $\partial^{-3}=\left(
                                                                                                     \begin{array}{c}
                                                                                                       f^{-2} \\
                                                                                                       -d^{-2} \\
                                                                                                     \end{array}
                                                                                                   \right)
$, $\partial^{-2}=\left(
                    \begin{array}{cc}
                      d^{-2} & f^{-1} \\
                      0 & -d^{-1} \\
                    \end{array}
                  \right)
$, and $\partial^{-1}=\left(d^{-1},f^{0}\right)$.
\end{lem}
Let $A^!$ be the quadratic dual of $A$. As graded vector spaces $A^!_0\cong \k$, $A^!_1\cong V^*$ and $A^!_2\cong \k r^*$, where $r^*\in (\k r)^*$ defined by $r^*(r)=1$. The multiplication on $A^!$ is given by: $\alpha\beta=(a_1,\dots,b_n)M(b_1,\dots,b_n)^tr^*$, for $\alpha=a_1x_1^*+\cdots+a_nx_n^*$ and $\beta=b_1x_1^*+\cdots+b_nx_n^*$ in $V^*$ (cf. \cite[Section 3]{HVZ2}), where $\{x^*_1,\dots,x^*_n\}$ is the basis of $V^*$ dual to the basis $\{x_1,\dots,x_n\}$. Write $E^i(B):=\Ext_B^i({}_B\k,{}_B\k)$ and $E(B):=\op_{i\ge0}E^i(B)$. Then $E(B)$ is a graded algebra with the degree $i$ component $E^i(B)$. The minimal projective resolution of ${}_B\k$ above implies that, as graded vector spaces,
\begin{equation}\label{eq3}
 E(B)\cong A^!\op A^!(-1).
\end{equation}
We write an element in $E(B)$ as $(\alpha,\beta)$ for some $\alpha, \beta\in A^!$, and demote the Yoneda product on $E(B)$ by $(\alpha,\beta)*(\alpha',\beta')$.

\begin{prop} \label{prop1} Assume $\delta(r)=0$. Then $A[z;\bar{\delta}]$ is a 3-CY algebra.
\end{prop}
\proof By \cite[Proposition 3.3]{HVZ1}) in the Koszul case, $B=A[z;\bar{\delta}]$ is Calabi-Yau if and only if $E(B)$ is a graded symmetric algebra. Recall that a finite dimensional graded algebra $E=\op_{i\ge0}E^i$ is graded symmetric if there is an integer $d$ and a homogeneous
nondegenerate bilinear form $\langle-,-\rangle:E\times
E\longrightarrow \k(d)$ such that $\langle \alpha\beta,\gamma\rangle=\langle
\alpha,\beta\gamma\rangle$ and $\langle \alpha,\beta\rangle=(-1)^{ij}\langle
\beta,\alpha\rangle$ for all homogeneous elements $\alpha\in E^i,\beta\in E^j$ and $\gamma\in E^k$. Since the global dimension of $B$ is 3 and $\dim E^3(B)=1$, $E(B)$ is graded symmetric if and only if, for all elements $\Phi\in E^1(B), \Theta\in E^2(B)$, $\Phi*\Theta=\Theta*\Phi$. Let $\Phi=(\alpha,k)$ with $\alpha\in A^!_1=V^*$ and $k\in\k$, and $\Theta=(r^*,\beta)$ with $\beta\in V^*$. The element $\Phi$ induces a $B$-module morphism $g:B\ot V\op B\longrightarrow {}_B\k$ by $g(1\ot x,1)=\alpha(x)+k$ for all $x\in V$, and the element $\Theta$ induces a $B$-module morphism $h:B\ot \k r\op B\ot V\longrightarrow {}_B\k$ by $h(1\ot r,1\ot x)=1+\beta(x)$ for all $x\in V$.
 Consider the following diagram:
$$\xymatrix{0\ar[r]&B\ot \k r\ar[r]^{\partial^{-3}}\ar[d]_{g_{2}} &B\ot \k r\op B\ot V\ar[d]_{g_{1}}\ar[r]^{\quad\partial^{-2}}&B\ot V\op B\ar[d]_{g_{0}}\ar[r]^{\quad\partial^{-1}}\ar[dr]^{g} &B\cdots\\
 \cdots\ar[r]&B\ot \k r\op B\ot V\ar[r]^{\partial^{-2}} & B\ot V\op B\ar[r]^{\partial^{-1}}& B\ar[r]&{}_B\k,
}$$
where the vertical arrows are $B$-module morphisms defined as follows. As before, we write $r=\sum_{i=1}^nu_i\ot x_i$ with all $u_i\in V$, and assume $\delta(x_i)=\sum_{j=1}^n y_{ij}\ot x_j$ for all $i=1,\dots,n$ with all $y_{ij}\in V$.
Then
{\small\begin{eqnarray}
\nonumber g_0(1\ot x_j,1)&=&\alpha(x_j)1+k1;\\
\nonumber g_1(1\ot r,1\ot x_j)&=&\left(\sum_{i=1}^n1\ot u_i\alpha(x_i)-1\ot k x_j-\sum_{i=1}^n1\ot y_{ji}\alpha(x_j),\ \alpha(x_j)1\right);\\
 \nonumber g_2(1\ot r)&=&(1\ot k r,\ \sum_{i=1}^n1\ot u_i\alpha (x_i)),
\end{eqnarray}} for all $j=1,\dots,n$. Since $\delta(r)=0$, it follows that $\displaystyle\sum_{i=1}^n\delta(u_i) \ot x_i+\sum_{i=1}^n\sum_{j=1}^n u_i\ot y_{ij}\ot x_j=0$. Applying the linear map $1\ot 1\ot \alpha$ to this identity, one obtains:
 \begin{equation}\label{eqq4} \sum_{i=1}^n\delta(u_i)\alpha(x_i)+\sum_{i=1}^n\sum_{j=1}^nu_i\ot y_{ij}\alpha(x_j)=0.
\end{equation}
Using Equation (\ref{eqq4}) and  the following computations:
{\small\begin{eqnarray}
 \nonumber g_1\circ\partial^{-3}(1\ot r) &=& g_1(z\ot r,-r) \\
 \nonumber  &=&\left( \sum_{i=1}^nz\ot u_i\alpha(x_i)+\sum_{i=1}^nu_i\ot k x_i+\sum_{i=1}^n\sum_{j=1}^nu_i\ot y_{ij}\alpha(x_j),-\sum_{i=1}^n u_i\alpha(x_i)\right),
\end{eqnarray}
\begin{eqnarray}
  \nonumber \partial^{-2}\circ g_2(1\ot r) &=& \partial^{-2} \left(1\ot k r,1\ot\sum_{i=1}^nu_i\alpha(x_i)\right) \\
  \nonumber &=& \left(k r+\sum_{i=1}^nz\ot u_i\alpha(x_i)-\sum_{i=1}^n\delta(u_i)\alpha(x_i),-\sum_{i=1}^n u_i\alpha(x_i)\right).
\end{eqnarray}}
we obtain the identity: $g_1\circ\partial^{-3}(1\ot r)=\partial^{-2}\circ g_2(1\ot r)$. Hence $g_1\circ\partial^{-3}=\partial^{-2}\circ g_2$. Similar computations show that the second square in the diagram commutes. The commutativity of the triangle in the diagram is obvious.
Thus, we have $h\circ g_2(1\ot r)=h(1\ot k r, \sum_{i=1}^n1\ot u_i\alpha(x_i))=k+\sum_{i=1}^n\beta(u_i)\alpha(x_i)$. By the definition of the Yoneda product,  we have $\Theta*\Phi=(r^*,\beta)*(\alpha,k)=kr^*+\beta\alpha$, where $\beta\alpha$ is the product in $A^!$. Similarly, we can show that $\Phi*\Theta=kr^*-\alpha\beta$. Now $A$ is Calabi-Yau, hence $A^!$ is graded symmetric; that is, $\alpha\beta=-\beta\alpha$ for all $\alpha,\beta\in A^!_1$. It follows that $\Phi*\Theta=\Theta*\Phi$. Therefore, $B=A[z;\bar{\delta}]$ is Calabi-Yau. \qed

The computation in the proof of Proposition \ref{prop1} has given us the formulas of the Yoneda product of $E(B)$.

\begin{cor} \label{cor1} As vector spaces, $E(B)\cong A^!\op A^!(-1)$. The Yoneda product of $E(B)$ is given as follows: for $\alpha,\beta\in A^!_1$ and $k,k'\in\k$, $$(r^*,\beta)*(\alpha,k)=(\alpha,k)*(r^*,\beta)=kr^*+\beta\alpha,$$ and $$(\beta,k')*(\alpha,k)=(\beta\alpha,k'\alpha-k\beta-(\beta\ot\alpha)\circ\delta),$$ where $r^*$ is the basis of $A^!_2$ such that $r^*(r)=1$.
\end{cor}
\proof The first identity is proved in the proof of Proposition \ref{prop1}. Keep the same notions as in the proof of Proposition \ref{prop1}. The element $(\beta,k')$ induces a $B$-module morphism $g':B\ot V\op B\longrightarrow {}_B\k$ by $g'(1\ot x,1)=\beta(x)+k'$ for all $x\in V$, and $(\beta\alpha,k'\alpha-k\beta-(\beta\ot\alpha)\circ\delta)$ induces a $B$-module morphism $f:B\ot \k r\op B\ot V\longrightarrow {}_B\k$ by $f(1\ot r,1\ot x_j)=\sum_{i=1}^n\beta(u_i)\alpha(x_i)+k'\alpha(x_j)-k\beta(x_j)-\sum_{i=1}^n\beta(y_{ji})\alpha(x_i)$ for all $j=1,\dots,n$. By the definition of Yoneda product, $(\beta,k')*(\alpha,k)$ is represented by $g'\circ g_1$. Now $g'\circ g_1(1\ot r,1\ot x_j)=\sum_{i=1}^n\beta(u_i)\alpha(x_i)-k\beta(x_j)-\sum_{i=1}^n\beta(y_{ji})\alpha(x_i)+k'\alpha(x_j)=f(1\ot r,1\ot x_j)$ for all $j=1,\dots,n$. Therefore the second identity follows. \qed

Let $\epsilon:A^!\to A^!$ be the automorphism of $A^!$ defined by $\epsilon(\alpha)=-\alpha$ for $\alpha\in A^!_1$ and $\epsilon(\beta)=\beta$ for all $\beta\in A^!_2$. Let ${}_\epsilon A^!$ be the graded $A^!$-bimodule whose right $A^!$-action is the regular action, and whose left $A^!$-action is twisted by the automorphism $\epsilon$; that is, for all $\gamma,\theta\in A^!$, the left $A^!$-action $\gamma\cdot \theta=\epsilon(\gamma)\theta$. Let $I={}_\epsilon A^!(-1)$, and let $E(A^!;I)$ be the trivial extension of $A^!$ by the $A^!$-bimodule $I$. By Corollary \ref{cor1}, $E(B)$ is isomorphic to $E(A^!;I)$.

\begin{cor} \label{cor2} The Yoneda algebra $E(B)$ is isomorphic to the trivial extension of $A^!$ by the $A^!$-bimodule $I$.
\end{cor}

\begin{exa} \label{exa1} {\rm Consider the Calabi-Yau algebra studied by Smith in \cite{Sm}. Let $\k \langle x_1,\dots,x_6\rangle$ be the free algebra generated by six elements. Let $A=\k \langle x_1,\dots,x_6\rangle/\langle r\rangle$, where $$r=(x_1,\dots,x_6)\left(
                                                                       \begin{array}{cccccc}
                                                                         0 & 0 & 0 & 0 & 0 & 1 \\
                                                                         0 & 0 & 0 & 0 & -1 & 0 \\
                                                                         0 & 0 & 0 & -1 & 0 & 0 \\
                                                                         0 & 0 & 1 & 0 & 0 & 0 \\
                                                                         0 & 1 & 0& 0 & 0 & 0 \\
                                                                         -1 & 0 & 0 & 0 & 0 & 0 \\
                                                                       \end{array}
                                                                     \right)
\left(
  \begin{array}{c}
    x_1 \\
    \vdots \\
    x_6 \\
  \end{array}
\right).$$ Define a derivation $\delta:\k \langle x_1,\dots,x_6\rangle\to \k \langle x_1,\dots,x_6\rangle$ by
$$\begin{array}{ccccccc}
  \delta(x_1)&=&x_4x_2-x_2x_4+x_3x_5-x_5x_3&&\delta(x_2)&=&x_1x_4-x_4x_1+x_3x_6-x_6x_3\\
  \delta(x_3)&=&x_5x_1-x_1x_5+x_6x_2-x_2x_6&&\delta(x_4)&=&x_2x_1-x_1x_2+x_5x_6-x_6x_5\\
  \delta(x_5)&=&x_1x_3-x_3x_1+x_6x_4-x_4x_6&&\delta(x_6)&=&x_2x_3-x_3x_2+x_4x_5-x_5x_4.
  \end{array}
$$ Then $\delta(r)=0$, and $B=A[z;\bar{\delta}]$ is 3-CY \cite{Sm}.}
\end{exa}

Keep the assumption that $\delta(r)=0$. Let $\widehat{V}=V\op \k z$. Then $B=A[z;\bar{\delta}]$ is a quotient algebra of $T(\widehat{V})$. Since $B$ is 3-CY, $B$ is defined by a superpotential \cite[Theorem 3.1]{Bo}. Let $\{x_1^*,\dots,x_n^*\}$ be the basis of $V^*$ dual to $\{x_1,\dots,x_n\}$. Recall that a {\it superpotential} is an element $w\in \widehat{V}\ot \widehat{V}\ot \widehat{V}$ such that $[\alpha w]=[w\alpha]$ for all $\alpha\in(\widehat{V})^*$, where $[\alpha w]=(\alpha\ot 1\ot 1)(w)$ and $[w\alpha]=(1\ot 1\ot \alpha)(w)$. Given a superpotential $w$, the {\it partial derivative} of $w$ by $x_i$ is defined by $\partial_{x_i}(w)=[x^*_iw]$ (cf. \cite{BSW}). By \cite[Theorem 3.1]{Bo}, there is a superpotential $w\in \widehat{V}\ot \widehat{V}\ot \widehat{V}$ such that $B\cong T(\widehat{V})/\langle \partial_{x_i}(w):i=0,\dots,x_n\rangle$ where $x_0=z$. We next show that the superpotential $w$ may be written out explicitly. For $i=1,\dots,n$, let $r_i=z\ot x_i-x_i\ot z-\delta(x_i)\in\widehat{V}\ot \widehat{V}$. Clearly $r,r_1,\dots,r_n$ are linearly independent in $\widehat{V}\ot \widehat{V}$, moreover $B\cong T(\widehat{V})/\langle r,r_1,\dots,r_n\rangle$. Before we construct the general form of the superpotentials, let us look at the following example.

\begin{exa}\label{exa2} {\rm Let $\k \langle x,y\rangle$ be the free algebra generated by two elements. Let $\delta:\k \langle x,y\rangle\to \k \langle x,y\rangle$ be a derivation defined by $\delta(x)=bx^2+cy^2$ and $\delta(y)=ax^2-bxy-byx$, where $(a,b,c)\in\k^3$. Let $r=xy-yx$. Then it is easy to see that $\delta(r)=0$. Therefore, $\delta$ induces a derivation $\bar{\delta}$ on $A=\k[x,y]$. Now $B=A[z;\bar{\delta}]$ is 3-CY. A straightforward verification shows that $w=yxz+zyx+xzy-xyz-zxy-yzx-ax^3+cy^3+bxyx+bx^2y+byx^2$ is a superpotential, and $B\cong \k \langle x,y,z \rangle/\langle\partial_x(w),\partial_{y}(w),\partial_{z}(w)\rangle$. Explicitly, the generating relations are $r_1=zy-yz-ax^2+byx+bxy, r_2=xx-zx+cy^2+bx^2$ and $r_3=yx-xy$.}
\end{exa}

\begin{prop} \label{prop3} Assume $\delta(r)=0$. Let $Q=\left(
                                   \begin{array}{ccc}
                                     -1 & 0  \\
                                     0 & M \\
                                   \end{array}
                                 \right)$, and let        $w=(z,x_1,\dots,x_n)Q\left(
                                                                          \begin{array}{c}
                                                                            r \\
                                                                            r_1 \\
                                                                            \vdots \\
                                                                            r_n
                                                                          \end{array}
                                                                        \right)$, where $M$ is an invertible $n\times n$ anti-symmetric matrix, and $r_i=z\ot x_i-x_i\ot z-\delta(x_i)\in \widehat{V}\ot \widehat{V}$ for all $i=1,\dots,n$. Then

{\rm (i)} $w$ is a superpotential;

{\rm (ii)} $A[z;\bar{\delta}]\cong T(\widehat{V})/\langle \partial_{x_i}(w):i=0,\dots,n\rangle$, where we set $x_0=z$.
 \end{prop}
\proof Let $\{m_{ij}|i,j=1,\dots,n\}$ be the entries of $M$. Then $r=\sum_{i,j=1}^nm_{ij}x_i\ot x_j$. Since $\delta(r)=0$, we have $\displaystyle\sum_{i,j=1}^nm_{ij}\delta(x_i)\ot x_j=-\sum_{i,j=1}^nm_{ij}x_i\ot \delta(x_j)$.  Let us compute the element $w$.
$$\begin{array}{ccl}
   w&=&-z\ot r+(x_1,\dots,x_n)M\left(
                                                                          \begin{array}{c}
                                                                          r_1\\
                                                                            \vdots \\
                                                                            r_n
                                                                          \end{array}
                                                                        \right)\\
   &=&\displaystyle-\sum_{i,j=1}^nm_{ij}z\ot x_i\ot x_j+\sum_{i,j=1}^nm_{ij}x_i\ot r_j\\
   &=&\displaystyle-\sum_{i,j=1}^nm_{ij}z\ot x_i\ot x_j+\sum_{i,j=1}^nm_{ij}x_i\ot z\ot x_j\\
   &&\displaystyle-\sum_{i,j=1}^nm_{ij}x_i\ot x_j\ot z -\sum_{i,j=1}^nm_{ij}x_i\ot \delta(x_j)\\
   &=&-\displaystyle\sum_{i,j=1}^nm_{ij}(z\ot x_i-x_i\ot z)\ot x_j-\sum_{i,j=1}^nm_{ij}x_i\ot x_j\ot z +\sum_{i,j=1}^nm_{ij}\delta(x_i)\ot x_j\\
   &=&-\displaystyle\sum_{i,j=1}^nm_{ij}(z\ot x_i-x_i\ot z-\delta(x_i))\ot x_j-\sum_{i,j=1}^nm_{ij}x_i\ot x_j\ot z\\
   &=&-r\ot z+(r_1,\dots,r_n)M^t\left(
                                                                          \begin{array}{c}
                                                                          x_1\\
                                                                            \vdots \\
                                                                            x_n
                                                                          \end{array}
                                                                        \right).
  \end{array}
$$
Now it is clear that $[x_i^*w]=[wx_i^*]$, and $\partial_{x_i}(w)=r_i$ for all $i=0,1,\dots,n$, where $r_0=r$. \qed

\section{Coherence of $A[z;\bar{\delta}]$}

Notation and notions are as in the previous section. By \cite[Theorem 0.2]{Z}, $A$ is Noetherian if and only if $\dim (V)=2$. Since $B$ is an Ore extension of $A$ in variable $z$, $B/Bz$ is isomorphic to $A$ as a graded left $B$-module. Since $A$ is not left Noetherian when $\dim (V)>2$, neither is $B$. Similarly, $B$ is not right Noetherian when $\dim(V)>2$. Summarizing the foregoing argument, we obtain the following property.

\begin{lem} $B=A[z;\bar{\delta}]$ is Noetherian if and only if $\dim (V)=2$.
\end{lem}

Piontkovski showed in \cite[Theorem 4.1]{Pi} that any connected graded algebra with a single quadratic relation is graded coherent. Hence $A$ is a graded coherent algebra. So, it is natural to ask whether $B$ is a graded coherent algebra. The answer is affirmative. However, the proof of this property is not trivial because an Ore extension of a coherent ring needs not be coherent. In fact, there is a commutative coherent ring $R$ such that the polynomial extension $R[z]$ is not coherent \cite{So}. Some other results about the coherence of polynomial rings may be found in \cite{GV}. Let us recall the definition of a graded coherent algebra.

A graded algebra $D$ is called a {\it graded left coherent} algebra if one of the following equivalent conditions is satisfied:
\begin{itemize}
  \item [(i)] every finitely generated graded left ideal of $D$ is finitely presented; that is, if $I$ is a graded left ideal of $D$ then there is a finitely graded free $D$-module $F$ and a surjective morphism $g:F\to I$ of graded modules such that $\ker g$ is also a finitely generated $D$-module;
  \item[(ii)] every finitely generated graded submodule of a finitely presented graded module is finitely presented;
  \item [(iii)] the category of all finitely presented graded left $D$-modules is an abelian category.
\end{itemize}
Similarly we can define a {\it graded right coherent} algebra. If a graded algebra is both graded left and right coherent, then it is called a {\it graded coherent} algebra.

Let $W=\op_{i\ge0}W_i$ be a graded vector space with $\dim (W_i)<\infty$ for all $i$. Recall that the Hilbert series of $W$ is defined to be the power series $H_W(t)=\sum_{i\ge0}\dim(W_i) t^i$.

\begin{lem}\label{lem2} Let $V$ be a vector space of dimension $n\ge4$ with basis $\{x_1,\dots,x_n\}$, and let \begin{equation}\label{eq4} M=\left(
\begin{array}{cccccc}
                                                                                                          0 & \cdots & 0 & 0 & \cdots & 1 \\
                                                                                                        \vdots &   & \vdots & \vdots &  & \vdots \\
                                                                                                          0  & \cdots  & 0  & 1 & \cdots & 0  \\
                                                                                                          0  & \cdots  & -1 & 0 & \cdots  & 0  \\
                                                                                                     \vdots  &   & \vdots  & \vdots &   & \vdots  \\
                                                                                                          -1 & \cdots  & 0  &0  & \cdots  & 0  \\
                                                                                                        \end{array}
                                                                                                      \right)\end{equation}
be the invertible $n\times n$ anti-symmetric matrix with entries in the anti-diagonal line $1$ or $-1$ and others 0.
Let $r=(x_1,\dots,x_n)M(x_1,\dots,x_n)^t$, and $A=T(V)/\langle r\rangle$. Let $\delta$ be a derivation on $T(V)$ of degree 1. We write $\delta(x_i)=\sum_{s,t=1}^nk^i_{st}x_s\ot x_t$ for all $i=1,\dots,n$. Assume that $k^i_{nn}=0$ for all $i=1,\dots,n$ and $\delta(r)=0$. Let $\bar{\delta}$ be the derivation on $A$ induced by $\delta$. Write $B=A[z;\bar{\delta}]$.
Then the following hold.
\begin{itemize}
  \item [(i)] Let $I$ be the ideal of $B$ generated by the elements $x_1,\dots,x_{n-1}$. Then $B/I\cong \k[X,Z]$, where $\k[X,Z]$ is the commutative polynomial algebra in variables $X$ and $Z$;
  \item [(ii)] Let $L=\k x_1\op\cdots\op\k x_{n-1}$ and $L'=\k x_2\op\cdots\op\k x_{n-1}$. Then, as left $B$-modules, $I\cong B\ot(L\op L'x_n\op L'x_n^2\op\cdots)$, where $L'x_n^k$ ($k\ge1$) is the vector space spanned by the elements $x_2x_n^k,\dots,x_{n-1}x_n^k$.
\end{itemize}
\end{lem}

\begin{con}\label{con} {\rm We call an $n\times n$ ($n\ge2$) invertible anti-symmetric matrix of the form (\ref{eq4}) a {\it standard anti-symmetric matrix}. If $M$ is an invertible anti-symmetric matrix, there is an invertible matrix $P$ such that $P^tMP$ is standard.}
\end{con}

\noindent{\it Proof of Lemma \ref{lem2}}. (i) By assumption, $\delta(x_n)=\sum_{s,t=1}^nk^n_{st}x_s\ot x_t$ and $k^n_{nn}=0$. Therefore $\bar{\delta}(x_n)\in I$ and $B/I$ is a commutative algebra.
There is an algebra morphism $g:k[X,Z]\longrightarrow B/I$ defined by $g(X)=x_n$ and $g(Z)=z$. Next, we want to construct an algebra morphism from $B/I$ to $\k[X,Z]$. As before, write $\widehat{V}=V\oplus\k z$. Firstly, we define $f:T(\widehat{V})\longrightarrow \k[X,Z]$ by letting $f(x_i)=0$ for all $i=1,\dots,n-1$, $f(x_n)=X$ and $f(z)=Z$. Denote by $\langle x_1,\dots,x_{n-1}\rangle$ and by $\langle z\ot x_n-x_n\ot z\rangle$  the ideals of $T(\widehat{V})$ respectively generated by $x_1,\dots,x_{n-1}$  and  by $z\ot x_n-x_n\ot z$. Obviously, $\langle x_1,\dots,x_{n-1}\rangle+\langle z\ot x_n-x_n\ot z\rangle\subseteq\ker f$. Recall that $B$ is a Koszul algebra and $B=T(\widehat{V})/J$ where $J=\langle r,z\ot x_1-x_1\ot z-\delta(x_1),\dots,z\ot x_n-x_n\ot z-\delta(x_n)\rangle$. Since $\delta(x_i)=\sum_{s,t=1}^nk^i_{st}x_s\ot x_t$ such that $k^i_{nn}=0$ for all $i=1,\dots,n$, it follows that $\delta(x_i)\in\langle x_1,\dots,x_{n-1}\rangle$ for all $i=1,\dots,n$. Hence $r,z\ot x_1-x_1\ot z-\delta(x_1),\dots,z\ot x_{n-1}-x_{n-1}\ot z-\delta(x_{n-1})\in \langle x_1,\dots,x_{n-1}\rangle\subseteq\ker f$. Now $z\ot x_n-x_n\ot z-\delta(x_n)\in \langle z\ot x_n-x_n\ot z\rangle+\langle x_1,\dots,x_{n-1}\rangle\subseteq\ker f$. Hence $J\subseteq\ker f$. Therefore, $f$ induces an algebra morphism $\overline{f}:B\longrightarrow \k[X,Z]$. Obviously, $\ker\overline{f}\supseteq I$. Hence $\overline{f}$ in turn induces an algebra morphism $\hat{f}:B/I\longrightarrow \k[X,Z]$. Now it is easy to see that $\hat{f}\circ g=id=g\circ\hat{f}$. The statement (i) follows.

(ii) Here we make use of the technique from \cite[Proposition 7.3]{Sm}. Let $\mu:B\ot B\to B$ be the multiplication of $B$. Then the restriction of $\mu$ defines a left $B$-module morphism (also denoted by $\mu$): $$\mu:B\ot(L\op L'x_n\op L'x_n^2\op\cdots)\longrightarrow I.$$ We claim that $\mu$ is surjective. In fact, if we can show that the image $I'=\text{im}(\mu)$ is also an ideal of $B$, then $I=I'$. So, it suffices to show that $I'x_n\subseteq I'$ and $I'z\subseteq I'$. Following the generating relation of $A$, we have $x_1x_n=x_nx_1+(x_2x_{n-1}-x_{n-1}x_2)+\cdots+(x_{\frac{n}{2}}x_{\frac{n}{2}+1}-x_{\frac{n}{2}+1}x_{\frac{n}{2}})\in BL\subseteq I'$. Therefore $I'x_n\subseteq I'$. In particular, $\bar{\delta}(x_i)\in I'$ for all $i=1,\dots,n$. On the other hand, since $x_iz=zx_i-\bar{\delta}(x_i)$, it follows that $x_iz\in I'$ for all $i=1,\dots,n-1$. For $2\leq i\leq n-1$, we have $x_ix_nz=x_i(zx_n-\bar{\delta}(x_n))=x_izx_n-x_i\bar{\delta}(x_n)\in I'x_n+x_iI'\subseteq I'$. Now assume $x_ix_n^jz\in I'$ for all $j<k$ and $2\leq i\leq n-1$. Then $$x_ix_n^kz=x_ix_n^{k-1}(zx_n-\bar{\delta}(x_n))=(x_ix_n^{k-1}z)x_n-x_ix_n^{k-1}\bar{\delta}(x_n)\in I'x_n+x_ix_n^{k-1}I'\subseteq I'.$$ Hence $I'z\subseteq I'$. The claim follows.
To show that $\mu$ is injective, we only need to compare the Hilbert series of $I$ and that of $F:=B\ot(L\op L'x_n\op L'x_n^2\op\cdots)$. Write $W=L\op L'x_n\op L'x_n^2\op\cdots$. Clearly $H_F(t)=H_B(t)\cdot H_W(t)$. We have $$H_W(t)=(n-1)t+(n-2)t^2+(n-2)t^3+\cdots=((n-1)t-t^2)(1-t)^{-1}.$$ The exact sequence $0\to I\to B\to B/I\to0$ implies $H_I(t)=H_B(t)-H_{B/I}(t)$. Since $B$ is Koszul of global dimension 3, it follows that $H_B(t)=\left(1-(n+1)t+(n+1)t^2-t^3\right)^{-1}$ by \cite[Theorem 5.9]{Sm1} and the isomorphism (\ref{eq3}) of the previous section. By (i), $H_{B/I}(t)=(1-t)^{-2}$. Hence
$$\begin{array}{ccl}
 H_{I}(t)&=&\left(1-(n+1)t+(n+1)t^2-t^3\right)^{-1}-(1-t)^{-2}\\
 &=&  \left(1-(n+1)t+(n+1)t^2-t^3\right)^{-1}\cdot\left((n-1)t-t^2\right)(1-t)^{-1}\\
&=&H_B(t)\cdot H_W(t)\\
&=&H_F(t).
 \end{array}
$$ Therefore $\mu$ is injective. So, (ii) follows. \qed

{\it\textbf{ Proof of the statement (iii) of Theorem \ref{thm}}}. If $n=2$, then $A=\k[x_1,x_2]$. We obtain that $B=A[z;\bar{\delta}]$ is Noetherian, and hence coherent. Now assume $n\ge4$. We only prove the statement when $j=n$ in the assumption, that is, $k^i_{nn}=0$ for all $i=1,\dots,n$. When $j\neq n$, the statement can be proved similarly. By Lemma \ref{lem2}, there is an exact sequence $0\longrightarrow I\longrightarrow B\longrightarrow B/I\longrightarrow0$ such that $B/I$ is a polynomial algebra in two variables and $I$ is a free graded left $B$-module. By \cite[Proposition 3.2]{Pi}, $B$ is graded right coherent. Note that the left version of Lemma \ref{lem2}(ii) holds too. Hence $B$ is also graded left coherent. \qed

As a special case of the statement (iii) of Theorem \ref{thm}, we have the following result, which can be viewed as a noncommutative version of \cite[Theorem 4.3]{GV}.

\begin{prop}\label{prop4} Let $A$ be a connected graded 2-CY algebra. Then $A[z]$ is a graded coherent algebra.
\end{prop}
\proof By \cite[Theorem 0.1]{Z} (also, cf. \cite[Proposition 3.4]{Ber}), $A$ is defined by an invertible anti-symmetric matrix $M$, that is, $A=T(V)/\langle r\rangle$ with $r=(x_1,\dots,x_n)M(x_1,\dots,x_n)^t$. For an invertible anti-symmetric matrix $M$, there is an invertible matrix $P$ such that $P^tMP$ is a standard invertible anti-symmetric matrix. Then the algebras defined by $M$ and $P^tMP$ respectively are isomorphic to each other. Hence we may assume that the anti-symmetric matrix $M$ itself is standard. Now by (iii) of Theorem \ref{thm}, we see that $A[z]$ is graded coherent. \qed

Now assume that $B=A[z;\bar{\delta}]$ is graded coherent. We may form a noncommutative projective space from $B$. Following \cite{Po}, we denote by $\text{coh}B$ the category of all finitely presented graded left $B$-modules, and by $\text{fdim}B$ the category of all finite dimensional graded left $B$-modules. Since $B$ is graded coherent, fdim$B$ is a Serre subcategory of coh$B$. Hence the quotient category $$\text{cohproj}B:=\text{coh}B/\text{fdim}B$$ is also an abelian category. Since $B$ is Koszul and 3-CY, $B$ is Artin-Schelter regular with Gorenstein parameter $-3$. Hence the Beilinson algebra of $B$ (for the terminology, see \cite[Definition 4.7]{MM}) is $$\nabla B=\left(
                          \begin{array}{ccc}
                            \k & B_1 & B_2 \\
                            0 & \k & B_1 \\
                            0 & 0 & \k \\
                          \end{array}
                        \right).
$$ Let $\text{mod}\nabla B$ be the category of finite dimensional left $\nabla B$-modules. Then by \cite[Theorem 4.14]{MM}, we have the following corollary.

\begin{cor} If the conditions of Theorem \ref{thm} are satisfied, then there is an equivalence of triangulated categories: $$D^b(\text{\rm cohproj}B)\cong D^b(\text{\rm mod}\nabla B),$$ where $D^b(-)$ is the bounded derived category of the corresponding abelian category.
\end{cor}

\vspace{5mm}

\subsection*{Acknowledgement}
The authors are very grateful to the referee for his/her careful reading of the manuscript, numerous comments and suggestions, which has greatly improved the paper. In particular, Remark \ref{rrem} was suggested to the authors by the referee. The work is supported by an FWO-grant and grants from NSFC (No. 11171067) and NSF of Zhejiang Province (No. LY12A01013). This work is also in part supported by SRF for ROCS, SEM.

\vspace{5mm}

\bibliography{}

\end{document}